\documentstyle[12pt]{article}

\newtheorem{thm}{Theorem}[section]

\newtheorem{coro}[thm]{Corollary}
\newtheorem{lema}[thm]{Lemma}
\newenvironment{pf}    {\par\noindent{\bf {Proof}\/. }\nopagebreak\normalsize}     {\hfill\linebreak[0]\hspace*{\fill}$\square$\\[1pt]}

\input amssym.def
\input amssym.tex

\begin{document}
\title{Exceptional Points in the \\ Elliptic-Hyperelliptic Locus}
\author{Ewa Tyszkowska \thanks{Institute of Mathematics, University of Gda\'nsk,  Wita Stwosza 57, 80-952 Gda\'nsk,  Poland.  E-mail address: ewa.tyszkowska@math.univ.gda.pl.  Supported by BW 5100-5-0289-7.} 
\\ and \\ Anthony Weaver
 \thanks{Department of Mathematics and Computer Science, Bronx Community College, City University of New York, University Avenue and West 181st Street, Bronx, New York, 10453, U.S.A.  E-mail address:  anthony.weaver@bcc.cuny.edu.  Supported by PSC CUNY Research Award 35-9370.}
}
\maketitle 
\begin{abstract} 
 An exceptional point in the moduli space of compact Riemann surfaces is a unique surface class whose full automorphism group acts with a triangular signature.  A surface admitting a conformal involution with quotient an elliptic curve is called elliptic-hyperelliptic; one admitting an anticonformal involution is called symmetric.  In this paper, we determine, up to topological conjugacy, the full group of conformal and anticonformal automorphisms of a symmetric exceptional point in the elliptic-hyperelliptic locus.  We determine the number of ovals of any symmetry of such a surface.  We show that while the elliptic-hyperelliptic locus can contain an arbitrarily large number of exceptional points, no more than four are symmetric. 
\end{abstract}
{\it Key words.} elliptic-hyperelliptic surface, symmetric surface, automorphism group.

\section{Introduction}

Exceptional points in  moduli space are unique surface classes whose full group of conformal automorphisms acts with a triangular signature.  Such points are of interest for many reasons, one being that their defining equations (as algebraic curves) have coefficients in a number field \cite{[Bel]}.  
Determining the  exceptional points in ${\cal M}_g$ (the moduli space in genus $g$) is not a simple matter.  Although there are just finitely many possible triangular signatures satisfying the Riemann-Hurwitz relation in genus $g$, not  all of them correspond to a group action.  Furthermore, several distinct groups may act with the same signature, or one group may have topologically distinct actions with the same signature.     The problem can be attacked piecemeal by restricting attention to certain subloci in ${\cal M}_g$.

The   
{\it $n$-hyperelliptic locus} ${\cal M}_g^n \subseteq {\cal M}_g$ consists of surfaces admitting a conformal involution  (the {\it $n$-hyperelliptic involution}),  with  quotient  a surface of genus $n$. 
When
$n=0$, these are the classical hyperelliptic surfaces.  When $n=1$,  these are the {\it elliptic-hyperelliptic surfaces}.  If $g > 4n+1$,  the $n$-hyperelliptic involution is unique and central in the full group of conformal automorphisms of the surface, and this puts a strong structural restriction on any larger group of automorphisms of the surface.    The $n$-hyperelliptic loci are 
important in building up a stratification of
 ${\cal M}_g$, since  if
$g>4n+3$, the intersection ${\cal M}_g^n \cap {\cal M}_g^{n+1}$ is empty  (\cite{FK92}, Cor V.1.9.2; see also \cite{Ri93}, \cite{We07}).

 For
$n>1$,  there is a constant bound $168(n-1)$ on the order of the full automorphism group of an $n$-hyperelliptic surface of genus $g$.  It follows that for sufficiently large $g$, and $n>1$, ${\cal M}_g^n$  contains no exceptional
points.  The number of  exceptional points  in ${\cal M}_g^0$ is
always between 
$3$ and $5$ (inclusive) and is precisely $3$ for 
all $g >30$ \cite{We04}.Ê By contrast, we shall show (Theorem~\ref{T:infsequence}) that  for infinitely many $g$, the number of exceptional points in ${\cal M}_g^1$ is larger than any pre-assigned positive integer (but, also, for infinitely many $g$, the number of exceptional points  in ${\cal M}_g^1$ is $0$).  ÊÊ 

A {\it symmetry} of a Riemann surface  is an antiholomorphic
involution; a surface is {\it symmetric} if it admits a symmetry.    Under the correspondence
between curves and surfaces, the fact that a surface $X$ is
symmetric means that the corresponding curve is definable  over
${\Bbb R}$.  In the group of conformal and anticonformal
automorphisms of $X$, non-conjugate symmetries correspond bijectively to real curves which are   non-isomorphic (over ${\Bbb R}$), and whose complexifications are birationally equivalent to $X$.    It is natural to ask which exceptional points are also symmetric.  In the $0$-hyperelliptic locus with $g>30$, the answer is: all of them.  In Section~\ref{S:exceptional}  we show that if the elliptic-hyperelliptic locus contains exceptional points, at most four of them are also symmetric.

If $X$ has genus $g$,  and $\varrho$ is a symmetry of $X$, 
 the set of fixed points  ${\rm Fix}(\varrho)$ of $\varrho$
consists of $k$ disjoint Jordan curves called {\it ovals}, where  $0 \leq k \leq g+1$,  by a theorem of Harnack
\cite{[Harn]}.  In Section~\ref{S:symmetries}, using the topological classification of conformal actions on elliptic-hyperelliptic
Riemann surfces given in \cite{[T]} (with a corrigendum supplied in Section~\ref{S:fullcorrigendum}), together with a result of Singerman characterizing  symmetric exceptional points \cite{[Sing1]},  and  
 a formula of Gromadzki \cite{[G1]}, we give  presentations of the full group of conformal and anticonformal automorphisms of  symmetric exceptional points in  ${\cal M}_g^1$, $g>5$, and  count   the number of ovals of each conjugacy class of  symmetry in such a group.

The outline of the paper is as follows.  In Section~\ref{S:NEC} we give necessary preliminaries on NEC groups.  In  Section~\ref{S:elliptic},  we study actions, reduced by the $1$-hyperelliptic involution, on the quotient elliptic curves; the algebra and number theory of the Gaussian and Eisenstein integers yield natural presentations of the reduced groups.   In  Section~\ref{S:triangsymm} we give the topological classification of conformal actions with triangular signature on symmetric elliptic-hyperelliptic surfaces.  In Section~\ref{S:fullcorrigendum} we state a corrigendum to a theorem in \cite{[T]}, which contained errors with respect to maximality of actions.  We thus obtain a complete and correct classification (Theorem~\ref{T:corrigendum}) of full conformal actions on elliptic-hyperelliptic surfaces.  In the final two sections,  we determine, up to topological conjugacy, the full group of conformal and anticonformal automorphisms of the  symmetric exceptional points in ${\cal M}_g^1$,  and count the ovals corresponding to the symmetries in such groups (Theorems~\ref{T:genus} and \ref{T:fullgroups}).

We  call attention to the related papers  \cite{BC97}, \cite{BCP04},  \cite{BS84}, \cite{[BS]}.

\section{NEC groups}\label{S:NEC}
Every compact Riemann surface $X$ of genus $g\geq 2$ can
be represented as the orbit space of the hyperbolic plane
$\mathcal{H}$ under the action of a discrete, torsion-free group $\Gamma$, called a {\it surface group of genus $g$},  consisting of orientation-preserving isometries of ${\mathcal H}$, and isomorphic to the fundamental group of $X$.   
   Any group of conformal and anticonformal automorphisms  of $X ={\mathcal H}
/\Gamma$ can be represented as $\Lambda /\Gamma$, where $\Lambda$ is a {\it non-euclidean crystallographic (NEC) group} 
 containing $\Gamma$ as a normal subgroup. An NEC  group
is a co-compact discrete subgroup of  the full group $\mathcal{G}$ of isometries (including those which reverse orientation)  of
$\mathcal{H}$.   Let $\mathcal{G}^+$ denote a subgroup of
$\mathcal{G}$ consisting of orientation-preserving isometries.
An NEC group is called a {\it Fuchsian group } if it is
contained in $\mathcal{G}^+$,  and a {\it proper NEC group}
otherwise. 

 Wilkie \cite{Wi66} and Macbeath \cite{Ma66} associated to every NEC group a
{\it signature} which determines its algebraic structure and the geometric nature of its action.  
It has the form
\begin{equation}\label{sig}
(g;\pm ;[m_1,\ldots,m_r];\{(n_{11},\ldots
,n_{1s_1}),\ldots ,(n_{k1},\ldots ,n_{ks_k})\}),
\end{equation}
where the numbers $m_i\geq 2$ are called the {\it proper periods}, the
brackets $()$ (which may be empty) are called the {\it period cycles}, the
numbers $n_{ij}\geq 2$ are called the {\it link periods}, and $g\geq 0$ is
 the {\it orbit genus}.   An NEC group with 
signature of the form $(g;[-];\{(-),\ldots ,(-)\})$ is called a {\it surface
NEC group of genus $g$}.   A Fuchsian group is an NEC group with
 signature of the form 
\begin{equation} (g;+;[m_1,\ldots ,m_r];\{-\}).
\end{equation}
In the particular case $g=0$ we shall write briefly $[m_1,\ldots
,m_r]$.   A group with signature $[m_1, m_2, m_3]$ is called a {\it triangle group}, and the signature is called ${\it triangular}$.   If $\Lambda$ is a proper NEC group with the signature
(\ref{sig}),  its {\it canonical Fuchsian subgroup} $\Lambda
^+=\Lambda \cap \mathcal{ G}^+$ has the signature
\begin{equation}\label{sig-fuch}
 \!\!\!\!\!(\gamma;+;[m_1,m_1,\ldots
,m_r,m_r,n_{11},\ldots n_{1s_1},\ldots ,n_{k1},\ldots
,n_{ks_k}];\{-\}),
\end{equation} 
where $\gamma= \alpha g+k-1$ and $\alpha =2$ if the sign is $+$ and $\alpha
=1$ otherwise. The group with the signature (\ref{sig}) has a
presentation given by generators:
$$\begin{tabular}{lll}
$($i$)$ & $x_i,i=1,\ldots ,r,$ & (elliptic generators)\\
$($ii$)$ & $c_{ij},i=1,\ldots,k;j=0,\ldots s_i,$& (reflection generators) \\
$($iii$)$ & $e_i,i=1,\ldots ,k, $ & (boundary generators)\\
$($iv$)$  & $a_i,b_i,i=1,\ldots g\;{\rm if \;the \;sign\; is}\; +,$& (hyperbolic generators)\\
  & $d_i,i=1,\ldots g\;{\rm if \;the \;sign\; is}\; -,$ & (glide reflection generators)\\
\end{tabular}
$$
and relations 
$$\begin{tabular}{ll}
$\!\!\!\!\!\!\!\!\!\!\!$$($1$)$ & $x_i^{m_i}=1,i=1,\ldots ,r,$ \\
$\!\!\!\!\!\!\!\!\!\!\!$$($2$)$ & $c_{is_i}=e_i^{-1}c_{i0}e_i,i=1,\ldots,k ,$\\
$\!\!\!\!\!\!\!\!\!\!\!$$($3$)$ & $c_{ij-1}^2=c_{ij}^2=(c_{ij-1}c_{ij})^{n_{ij}}=1,i=1,\ldots ,k;j=1,\ldots ,s_i,$\\
$\!\!\!\!\!\!\!\!\!\!\!$$($4$)$ & $x_1\ldots x_re_1\ldots
e_ka_1b_1a_1^{-1}b_1^{-1}\ldots
a_gb_ga_g^{-1}b_g^{-1}=1,$\\
   & $x_1\ldots x_re_1\ldots e_kd_i^2\ldots d_g^2=1.$\\
\end{tabular}
$$

Any system of generators of an NEC group satisfying the above
relations will be called {\it a canonical system} of generators.

Every NEC group has a fundamental region, whose hyperbolic area
is given by \begin{equation}\label{area} \mu
(\Lambda)=2\pi(\alpha g+k-2+\sum _{i=1}^r(1-1/m_i)+1/2\sum
_{i=1}^k\sum _{i=1}^{s_i}(1-1/n_{ij})),
\end{equation} where $\alpha$ is defined as in (\ref{sig-fuch}).
It is known that an abstract group with the presentation given by
the generators $(i)-(iv)$ and the relations $(1)-(4)$ can be
realized as an NEC group with the signature (\ref{sig}) if and
only if the right-hand side of (\ref{area}) is positive. If
$\Gamma$ is a subgroup of finite index in an NEC group $\Lambda$
then it is an NEC group itself and the {\it Riemann-Hurwitz relation} is
\begin{equation}\label{E:RiemHur}
[\Lambda :\Gamma]=\mu (\Gamma)/\mu (\Lambda).
\end{equation}

\section{Conformal actions on elliptic curves}\label{S:elliptic}

Every elliptic curve is a quotient of the  additive group  of
${\Bbb C}$ by a {\it lattice} ${\cal L} = {\cal L}(\tau_1,
\tau_2) \subset {\Bbb C}$ of co-finite area with basis $\{\tau_1,
\tau_2\} \subset {\Bbb C}$, and {\it modulus\/} $\mu =\tau_1/
\tau_2 \notin {\Bbb R}$  (see, e.g., \cite{JS}, \S 5.8).   ${\Bbb
C}/{\cal L}$ and ${\Bbb C}/{\cal L}^\prime$ are conformally
equivalent if and only ifÊ the moduli $\mu$ and $\mu^\prime$ of
${\cal L}$ and ${\cal L}{^\prime}$ are equivalent under a
M\"obius transformation defined by an integer matrixÊÊ of
determinant $1$.ÊÊ In such a case, the lattices ${\cal L}, {\cal
L}^\prime$ are {\it similar\/},  which means ${\cal L}^\prime =
\lambda{\cal L}$ for some $\lambda \in {\Bbb
C}\setminus\{0\}$.ÊÊ  If ${\cal L}$ is {\it self-similar} ($
\lambda{\cal L}={\cal L}$, $\lambda \ne 1$), then ${\Bbb C}/{\cal
L}$ admits a conformal automorphism with a fixed point.ÊÊ  Ê
Every lattice ${\cal L}={\cal L}(\tau_1, \tau_2)$ isÊÊ
self-similar by ${\cal L}= -1{\cal L}$,  thusÊÊ every elliptic
curve admits the conformal involution $ h_{ -1}:ÊÊ[z] \mapsto
[-z] $, whereÊÊ $[z]$ denotes the coset $z + {\cal L}$.ÊÊ    If
${\cal L}$ has modulus $i$,ÊÊ  ${\cal L}=i{\cal L}$ and $h_{-1}$
is the {\it square} of the automorphism $ h_i:\quad [z] \mapsto
[iz]$ ÊÊof order $4$.ÊÊ  Similarly, if ${\cal L}$ has modulus
$\omega = e^{i\pi/3}$, ${\cal L}=\omega{\cal L}$ (since
$\omega^2=\omega-1$) and $h_{-1}$ is the {\it cube} of the
automorphism $h_{\omega}:\quad [z] \mapsto [\omega z] $ of order
$6$.ÊÊ   Conjugates of $h_i$, $h_{\omega}$ and their powers
areÊÊ  the only automorphisms with fixed points on elliptic
curves.  (This fact is also known as the {\it crystallographic
restriction}, see \cite{Cox73}, \S 4.32.)

Lattices with moduli equivalent to $i$ or
$\omega$ are closed under  multiplication  and hence admit  ring structures.  
The rings ${\cal G} \simeq 
\{a+bi \mid a,b, \in {\Bbb Z}\}$ and ${\cal 
E}\simeq 
\{a+b\omega 
\mid a,b, \in {\Bbb Z}\},$ 
are known as the Gaussian and Eisenstein integers,
respectively (see, e.g., \cite{Cl84},  \S\S 174, 175).  Each ring has a 
multiplicative norm 
$\delta$ on the non-zero elements (see table below)
taking values  in the positive integers; the elements of
norm 
$1$ are calledÊÊ {\it units}.   Every non-zero, 
non-unit factorizes uniquely (up to reordering 
and unit multiples) into primes; primes which 
differ by a unit multiple are called {\it 
associates}.ÊÊ    Every ideal is {\it 
principal}; a {\it 
prime ideal\/} is one which, when it contains a 
product, also contains one of the factors.ÊÊ An 
element is a {\em prime\/} if and only 
if the principal ideal it generates is a prime 
ideal.ÊÊ  Associate primes generate the same prime 
ideal.ÊÊ Finally, every proper ideal factorizes 
uniquely as a product of proper prime ideals.ÊÊ A rational primeÊ $p$Ê {\it splits} in ${\cal G}$ or 
${\cal E}$ if it factorizes non-trivially.ÊÊ  
$p$ splits in ${\cal G}$ if and only if $p=2$
or $p \equiv 1$ mod $4$, for then $p$ is expressible in
the form 
$p=a^2 + b^2$, 
$a,b \in {\Bbb Z}$, and thus 
$p=(a+bi)(a-bi)$.ÊÊ  Similarly, a prime $p$ splits in
${\cal E}$ if and only if $p=3$ or $p \equiv 1$
mod $6$, for then $p$ is expressible in the form 
$p=a^2 + ab +b^2$, and thus 
$p=(a+b\omega)(a-b\omega^2)$.ÊÊ $2$, $3$ are the unique
splitting primes in ${\cal G}$, ${\cal E}$, respectively,
whose prime factors are associates.  See 
Table~\ref{Ta:rings}, below.
\begin{table}[h] 
\caption{The rings ${\cal G}$ and ${\cal
E}$}
\begin{tabular}{|c| c | c| c|} 
\hline
ring & units & norm& splitting primes \\ 
\hline 
${\cal G}$ & $\pm 1, \pm i$ & 
$\delta_{\cal G}(a+bi)=a^2 + b^2$ & $2,\quad p 
\equiv 1$ 
mod $4$\\ 
${\cal E}$ & $\pm 1, \pm \omega, \pm 
\omega^2$& $\delta_{\cal E}(a+b\omega)=a^2 + ab + 
b^2$& 
$3,\quad p 
\equiv 1$ mod $6$ 
\\ 
\hline 
\end{tabular}
\label{Ta:rings} 
\end{table}

If ${\mathcal L}' \subset {\Bbb C}$ is a lattice and ${\mathcal L}\subset {\mathcal L'}$ is a sublattice of finite
index,  the group ${\cal 
L}^\prime/{\cal L}=T$, abelian of rank $\leq
2$, acts by {\it translations}Ê  Êon 
${\Bbb C}/{\cal L}$ (that is, without fixed 
points or short orbits).ÊÊ 
Conjugates of translations are translations, 
hence a finite group of automorphisms of an elliptic curve   
has semidirect product structure 
$T \rtimes H$, where  $H$ is a cyclic group
generated by a power of $h_i$ or $h_\omega$.  If $|H|>2$, the moduli of ${\cal L}$ and ${\cal L}^\prime$
are both $i$ or both $\omega$ and ${\cal L}$ must be an {\it ideal} in
${\cal L}^\prime$:  otherwise the action of 
$H$ on ${\Bbb C}/{\cal L}^\prime$ cannot lift 
through the covering ${\Bbb 
C}/{\cal L}\rightarrow {\Bbb C}/{\cal 
L^\prime}$ (with Galois group $T$) to an action of $T \rtimes H$ on 
${\Bbb C}/{\cal L}$.   

 Let $ \mathcal{R}= \mathcal{G}$
or $ \mathcal{R}= \mathcal{E}$ and let $z \mathcal{R}$ denote the
principal ideal generated by $z \in \mathcal{R}$.  Let $\tilde
w$ denote the element $w+z \mathcal{R}$ in the additive group
${\mathcal R}/z {\mathcal R} = T$.  Since $T$ decomposes as a direct sum of abelian $p$-groups, the following theorem (communicated to us by Ravi S. Kulkarni)  leads to  a complete classification of finite groups of automorphisms of elliptic curves.

\begin{thm}\label{T:pSylow}
(a) If $p$ is a rational prime, $${{\cal G} \over p^k{\cal G}} =
\langle \tilde 1, \tilde{i}\rangle \simeq {\Bbb Z}_{p^k} \oplus
{\Bbb Z}_{p^k}; \qquad {{\cal E} \over p^k{\cal E}}=\langle
\tilde 1, \widetilde{ -\omega} \rangle \simeq  {\Bbb Z}_{p^k}
\oplus {\Bbb Z}_{p^k}.$$

(b) If $\pi \in {\cal R}\setminus {\Bbb Z}$ is a factor of a
splitting prime $p$,
\begin{enumerate}
\item $p>3:\:\: {{\cal R} \over \pi^k{\cal R}}= \langle \tilde 1
\rangle \simeq {\Bbb Z}_{p^k}$

\item $p=2:\:\: {{\cal G} \over \pi^{2k}{\cal G}}=\langle \tilde 1, \tilde i
\rangle \simeq {\Bbb Z}_{2^k} \oplus {\Bbb Z}_{2^k};\;\; {{\cal G}
\over \pi^{2k-1}{\cal G}}=\langle \widetilde{ 1-i}, \tilde i
\rangle \simeq {\Bbb Z}_{2^{k-1}} \oplus {\Bbb Z}_{2^k}$

\item $p=3:\:\: {{\cal E} \over \pi^{2k}{\cal E}}=\langle \tilde 1,
\widetilde{-\omega} \rangle \simeq {\Bbb Z}_{3^{k}} \oplus {\Bbb
Z}_{3^k};\;\; {{\cal E} \over \pi^{2k-1}{\cal E}}=\langle
\widetilde{1+\omega}, \tilde{\omega} \rangle \simeq {\Bbb
Z}_{3^{k-1}} \oplus {\Bbb Z}_{3^k}. $\end{enumerate}
\end{thm}

The most general finite automorphism group $\tilde G$ of an elliptic curve has the structure
\begin{equation}\label{E:ellipticaut}
{\Bbb Z}_n \oplus {\Bbb Z}_m \  \rtimes \  {\Bbb Z}_t, \quad m | n, \quad t | 4 \hbox{\ or \ } t | 6.
\end{equation}
If  $n/m \leq 3$, a presentation of $\tilde G$ is obtained by
arguments of the following type. Let $\tilde G_A=\langle\tilde
1,\tilde i\rangle\rtimes \langle h_i\rangle$, $\tilde
G_B=\langle\widetilde{1-i},\tilde i\rangle\rtimes \langle
h_i\rangle$, $\tilde G_C=\langle\tilde 1,\tilde \omega
\rangle\rtimes \langle h_\omega\rangle$, $\tilde
G_D=\langle\widetilde{1+\omega},\tilde \omega\rangle\rtimes
\langle h_\omega\rangle$, $\tilde G_E=\langle\tilde 1,\tilde
\omega \rangle\rtimes \langle h_{\omega^2}\rangle$ and $\tilde
G_F=\langle\widetilde {1+\omega},\tilde \omega\rangle\rtimes
\langle h_{\omega^2}\rangle$. Express $h_i$ as conjugation by an
element $c$,  and denote the elements $\tilde i,\tilde 1$ by $ x,
y$.  Then the relations $h_i(\tilde 1)=\tilde i,\;h_i(\tilde
i)=\tilde i^2=-\tilde 1$ become $ c y c^{-1}= x, c x c^{-1}=
y^{-1}$. In a similar way we obtain all the presentations given
in Table~\ref{Ta:reducedgroups}.

{\small \begin{table}[htdp] \caption{Some finite automorphism
groups of elliptic curves }
\begin{center}
\begin{tabular}[htp]{|l|l|} \hline
&\\[-3mm]
Group & Relators \\ \hline
$\tilde G_A=({\Bbb Z}_n \oplus {\Bbb Z}_n) \rtimes  {\Bbb Z}_4$ & $x^n, c^4,[x,y], cyc^{-1}x^{-1},cxc^{-1}y$\\[1mm]
$\tilde G_B=({\Bbb Z}_n \oplus {\Bbb Z}_{n/2})\rtimes {\Bbb Z}_4$ & $c^4,v^n,w^{n/2},
cvc^{-1}vw,cwc^{-1}v^{-2}w^{-1},
[v,w]$\\[1mm]
$\tilde G_C=({\Bbb Z}_n \oplus {\Bbb Z}_n)\rtimes {\Bbb Z}_6$ &
$y^n,c^6,cxc^{-1}y^{-1},cyc^{-1}y^{-1}x,[x,y]$\\[1mm]
$\tilde G_D=({\Bbb Z}_n \oplus {\Bbb Z}_{n/3})\rtimes {\Bbb Z}_6$ & $c^6,v^n,w^{n/3},
cvc^{-1}wv^{-2},cwc^{-1}wv^{-3},
[v,w]$\\[1mm]
$\tilde G_E=({\Bbb Z}_n \oplus {\Bbb Z}_n)\rtimes {\Bbb Z}_3$ & $c^3, x^n,cxc^{-1}y^{-1}x,cyc^{-1}x,[x,y]$\\[1mm]
$\tilde G_F=({\Bbb Z}_n \oplus {\Bbb Z}_{n/3})\rtimes {\Bbb Z}_3$ & $c^3,v^n,w^{n/3},
cvc^{-1}wv^{-1},cwc^{-1}w^2v^{-3}, [v,w]$\\[1mm]
 \hline
\end{tabular}
\end{center}
\label{Ta:reducedgroups}
\end{table}}

Groups of the more general type (\ref{E:ellipticaut}) ($t>2$) are obtained by adjoining a relation of the form $x^m=y^{mk}$ to one of the presentations in Table~\ref{Ta:reducedgroups}, where $k$ is a  positive integer such that $(k,n/m)=1$.  A finite group is thereby defined if and only if $k$ satisfies certain other conditions:  if $t=4$, we must have $k^2+1 \equiv 0$ mod ${n/m}$; if $t=3$ or $6$, we must have $k^2 - k + 1 \equiv 0$ mod ${n/m}$.  The existence of such $k$'s is equivalent to the number-theoretical condition that  all prime factors of $n/m$ split in the appropriate ring (${\cal G}$ if $t=4$, ${\cal E}$ if $t=3,6$).  

We note that  $2$ splits in ${\cal G}$ and $k^2 + 1 \equiv 0$ mod $2$ has the solution $k=1$, so we obtain a group isomorphic to $\tilde G_B$ by adjoining the relation $ x^{n/2}=
y^{n/2}$  to the presentation of $\tilde G_A$.  Similarly, $3$ splits in ${\cal E}$ and  $k^2 - k + 1 \equiv 0$ mod $3$ has the solution $k=2$, so we obtain a group isomorphic 
$\tilde G_D$ or
$\tilde G_F$ by adjoining the relation $x^{n/3}=y^{2n/3}$ to  the presentation of $\tilde G_C$ or $\tilde
G_E$, respectively.

\section{Triangular symmetric actions}\label{S:triangsymm}

If a group $G$ of conformal automorphisms of an
elliptic-hyperelliptic surface $X={\mathcal H}/\Gamma$ of genus
$g>5$ acts with  a triangular signature, then $G$ is isomorphic
to $\Lambda/\Gamma$, where $\Lambda$ has one of the fifteen
distinct triangular signatures
\begin{equation}\label{E:ehypsignatures}
[k,l,m]=[2\varepsilon_1,4\varepsilon_2,4\varepsilon_3],[3\varepsilon_1,3\varepsilon_2,3\varepsilon_3]
\hbox{\ or\ } [2\varepsilon_1,3\varepsilon_2,6\varepsilon_3], \quad
\varepsilon_i\in\{1,2\}
\end{equation}
where at least one  $\varepsilon_i$ is equal to $2$.  If the
signature of $\Lambda$ arises from $[2,4,4]$ (resp. $[2,3,6]$,
$[3,3,3]$) in this way, we shall say that $G$ determines a {\it
$(2,4,4)$-action} (resp., a {\it  $(2,3,6)$-,
$(3,3,3)$-action}).   Let $\theta:\Lambda \rightarrow G$ be an
epimorphism with kernel $\Gamma$. Then $\theta$ preserves the
(finite) orders of the  generators of $\Lambda$ and so $G$ is
generated by two elements $g_1$ and $g_2$ of orders $k$ and $l$
respectively whose product has order $m$. Singerman \cite{[Sing1]}
showed that $X$ is symmetric if and only if either of the maps
\begin {equation} \label{automorphism}
  g_1 \mapsto  g_1^{-1},  g_2 \mapsto  g_2^{-1}\quad {\rm or }\quad  g_1 \mapsto
g_2^{-1},  g_2 \mapsto  g_1^{-1}
 \end {equation}
induces an automorphism of $G$. Here we list all actions with
triangular signature on elliptic-hyperelliptic Riemann surfaces 
for which the group $G = \langle g_1, g_2 \rangle$ satisfies
condition (\ref{automorphism}). We call them {\it triangular
symmetric actions}.

\begin{thm}\label{T:action244} The topological type of a triangular symmetric $(2,4,4)$-action on an 
elliptic-hyper\-ellip\-tic Riemann surface is determined by a
finite group $G=A_\varepsilon^n$ or $G=B_{\varepsilon \delta}^n$
with presentation
\begin{equation}\label{one}
 \!\!\!\!\!A_\varepsilon^n=\langle x,y,c,\rho:
\rho^2, x^n\rho ^\varepsilon,c^4\rho^{\mu},[x,y]\rho^{\gamma},
cyc^{-1}x^{-1},cxc^{-1}y\rho^{\alpha},R\rangle,\end{equation}
 or
\begin{equation}\label{two}
{\bf \!\!\!\!\!\!\begin{array}{ll}B_{\varepsilon \delta}^n=\langle
w,v,c,\rho: \rho^2, c^4\rho^{\mu},v^n\rho^{\varepsilon},w^{n/2}\rho^\delta, cvc^{-1}wv\rho^{\alpha},& \\
\;\;\;\;\;\;\;\;\;\;\;\;\;\;\;\;\;\;\;\;\;\;\;\;\;\;\;\;\;\;\;\;\;\;\;\;\;\;cwc^{-1}v^{-2}w^{-1}\rho^{\alpha+\gamma},
[v,w]\rho^{\gamma},R\rangle, & \end{array}}\end{equation}
 a
Fuchsian group $\Lambda=\Lambda_{\alpha, \gamma, \mu}$  with
signature $[2(|\alpha-\mu|+1), 4(\mu + 1), 4(|\gamma - \mu| + 1)]$
generated by $x_1, x_2, x_3$,  and an epimorphism $ \theta:
\Lambda \rightarrow G$ defined by $\theta(x_1) = c^{-2}x$,
$\theta(x_2)=c$, $\theta(x_3) = y^{-1}c$ or
$\theta(x_1)=c^{-2}v,\theta(x_2)=c,\theta(x_3)=v^{-1}w^{-1}c$,
respectively, where $R$ is the set of relations making $\rho$
central and $\varepsilon,\delta\in\{0,1\}$. All non-equivalent
actions are listed in Table~\ref{Ta:244actions}. {\small
\begin{table}[htdp]
\caption{$(2,4,4)$-triangular symmetric actions}
\begin{center}
\begin{tabular}[!htp]{| c | c | c | c | c | c | c |} \hline
Case in \cite{[T]} & $\alpha$ & $\gamma$ & $\mu$& $n\equiv 0$ mod $4$ & $n\equiv 2$ mod $4$ & $n\equiv 1$  mod $2$\\
\hline
&&&&&&\\[-3mm]
$4.2$ & $0$ & $0$ & $1$ & $A_0^n,A_1^n,B_{00}^n,B_{01}^n$ & $A_0^n,A_1^n,B_{00}^n,B_{01}^n$     &$A_0^n,A_1^n$\\[1mm]
$4.3$ & $1$ & $0$ & $0$ & $A_0^n,A_1^n,B_{00}^n,B_{01}^n$ & $A_0^n,A_1^n,B_{10}^,B_{11}^n$     & \\[1mm]
$4.4$ & $0$ & $1$ & $0$ &$A_0^n,A_1^n,B_{00}^n,B_{01}^n$  &             &\\[1mm]
$4.5$ & $1$ & $1$ & $0$ & $A_0^n,A_1^n, B_{00}^n,B_{01}^n$ &            &\\[1mm]
$4.6$ & $1$ & $0$ & $1$ & $A_0^n,A_1^n,B_{00}^n,B_{01}^n$ & $A_0^n,A_1^n,B_{10}^n,B_{11}^n$     & \\[1mm]
\hline
\end{tabular}
\end{center}
\label{Ta:244actions}
\end{table}}
\end{thm}

\begin{thm}\label{T:action236}
The topological type of the triangular symmetric $(2,3,6)$-action
on an elliptic-hyperelliptic Riemann surface is determined by a
finite group of automorphisms $G=C_\varepsilon^n$ or
$G=D_{\varepsilon\varepsilon'}^n$ with the presentation
\begin{equation}\label{C236}
C_\varepsilon^n=\langle x,y,c,
\rho:\rho^2, y^n\delta^\varepsilon,c^6\rho^{\mu},cxc^{-1}y^{-1}\rho^{\alpha},
cyc^{-1}y^{-1}x\rho^{\beta},[x,y]\rho^{\alpha +\mu
},R\rangle
\end{equation}
 or
 \begin{equation}\label{D236}
  {\bf
\!\!\!\!\!\!\!\!\!\!\begin{array}{ll}D_{\varepsilon \varepsilon '
}^n=\langle
w,v,c,\rho:\rho^2, c^6\rho^{\mu},v^n\rho^{\varepsilon},w^{n/3}\rho^{\varepsilon
' },cvc^{-1}wv^{-2}\rho^{\alpha +\mu +\beta},
& \\
\;\;\;\;\;\;\;\;\;\;\;\;\;\;\;\;\;\;\;\;\;\;\;\;\;\;\;\;\;\;\;\;cwc^{-1}wv^{-3}\rho^{\beta+\mu},
[v,w]\rho^{\alpha+\mu},R\rangle , & \end{array}}\end{equation}
 a
Fuchsian group $\Lambda =\Lambda_{\alpha,\beta,\mu}$ with the
signature $[2(\alpha+1),3(\beta+1),6(\mu+1)]$ generated by
$x_1,x_2,x_3$,  and an epimorphism $\theta : \Lambda \to G$ defined
by $\theta (x_1)=c^{3}x,  \theta (x_2)=c^2y, \theta (x_3)=c$ or
$\theta (x_1)=c^3wv^{-1},\theta(x_2)=c^2v,\theta(x_3)=c$, 
respectively, where $R$ is the set of relations making $\rho$
central, $n$ is a positive integer and $\varepsilon\in\{0,1\}.$
All non-equivalent actions are listed in
Table~\ref{Ta:236actions}. {\small
\begin{table}[htdp]
\caption{$(2,3,6)$-triangular symmetric actions}
\begin{center}
\begin{tabular}{|c|c|c| c||c| c | c | c | c | c| c|}
\hline&&&& \multicolumn{7}{ |c |}{ $n$ mod $12$}
\\ \cline{5-11}
Case in \cite{[T]}& $\alpha$ & $\beta$ & $\mu$& $0$ & $ \pm 1$ & $\pm 2$ &
$\pm 3$ & $\pm 4$ & $\pm5$ & $ 6$\\ \hline
$6.2$ & $0$ & $1$ &  $0$ &$C_0^n,D_{00}^n$ &$C_1^n$ &$C_0^n$  &$C_1^n$, $D_{10}^n$  &$C_0^n$ &$C_1^n$ &$C_0^n,D_{00}^n$ \\[1mm]
$6.3$ & $0$ & $0$ &  $1$ &$C_0^n,D_{00}^n$ &      &$C_1^n$   &             &$C_0^n$ &      & $C_1^n,D_{11}^n$\\[1mm]
$6.4$ & $0$ & $1$ &  $1$ &$C_0^n,D_{00}^n$ &      &$C_1^n$   &             &$C_0^n$ &      &$C_1^n,D_{11}^n$  \\[1mm]
$6.5$ & $1$ & $0$ &  $0$ &$C_0^n,D_{00}^n$ &      &$C_1^n$   &             &$C_0^n$ &      &$C_1^n,D_{11}^n$ \\[1mm]
$6.6$ & $1$ & $1$ &  $0$ &$C_0^n,D_{00}^n$ &      &$C_1^n$   &             &$C_0^n$ &      & $C_1^n,D_{11}^n$ \\[1mm]
$6.7$ & $1$ & $0$ &  $1$ &$C_0^n,D_{00}^n$ &$C_1^n$  &$C_0^n$   & $C_1^n$, $D_{11}^n$  &$C_0^n$ &$C_1^n$  &$C_0^n,D_{00}^n$ \\[1mm]
$6.8$ & $1$ & $1$ &  $1$ &$C_0^n,D_{00}^n$ &$C_0^n$  &$C_0^n$   & $C_0^n$,$D_{01}^n$   &$C_0^n$ &$C_0^n$  &$C_0^n,D_{00}^n$ \\[1mm]
\hline
\end{tabular}
\end{center}
\label{Ta:236actions}
\end{table}}
\end{thm}

\begin{thm}\label{T:action333}
The topological type of the triangular symmetric $(3,3,3)$-action
on elliptic-hyperelliptic Riemann surface is determined by a
finite group of automorphisms $G=E_{\varepsilon}^n$ or
$G=F_{\varepsilon \varepsilon'}^n$ with the presentation
\begin{equation}\label{more333}E_{\varepsilon}^n=\langle x,y,c,\rho
:\rho^2, c^3\rho ^{\mu},
x^n\rho^\varepsilon,cxc^{-1}y^{-1}x,cyc^{-1}x\rho^{\beta},[x,y]\rho^{\gamma},R\rangle\end{equation}
or \begin{equation}\label{less333} {\bf
\!\!\!\!\!\!\!\!\!\!\begin{array}{ll}F_{\varepsilon \varepsilon '
}^n=\langle
w,v,c,\rho:\rho^2, c^3\rho^{\mu},v^n\rho^{\varepsilon},w^{n/3}\rho^{\varepsilon '}, cvc^{-1}wv^{-1}\rho^{\beta},& \\
\;\;\;\;\;\;\;\;\;\;\;\;\;\;\;\;\;\;\;\;\;\;\;\;\;\;\;\;cwc^{-1}w^2v^{-3}\rho^{\beta},
[v,w]\rho^{\gamma},R\rangle , & \end{array}}\end{equation} a
Fuchsian group $\Lambda_{\beta ,\gamma,\mu }$ with the signature
$[(\mu+1)3,(\beta+1)3,(|\gamma -\beta -\mu|+1)3]$ generated by
$x_1,x_2,x_3$ and an epimorphism $\theta:\Lambda \rightarrow G$
defined by
$\theta(x_1)=c,\theta(x_2)=c^{-2}x,\theta(x_3)=x^{-1}c$ or
$\theta(x_1)=c,\theta(x_2)=c^{-2}wv^{-1},\theta(x_3)=vw^{-1}c$,
respectively, where $R$ is the set of relations making $\rho$
central, $n$ is a positive integer and $\varepsilon,\varepsilon'
\in\{0,1\}$. All non-equivalent actions are listed in
Table~\ref{Ta:333actions}. {\small
\begin{table}[htdp]
\caption{$(3,3,3)$-triangular symmetric actions}
\begin{center}
\begin{tabular}{|c|c|c| c||c| c | c | c | c | c| c|}
\hline&&&& \multicolumn{7}{ |c |}{ $n$ mod $12$}
\\ \cline{5-11}
Case in \cite{[T]}& $\beta$ & $\gamma$ & $\mu$& $0$ & $ \pm 1$ & $\pm 2$ &
$\pm 3$ & $\pm 4$ & $\pm5$ & $ 6$\\ \hline
$3.2$ & $0$ & $1$ &  $0$ &$E_0^n,F_{00}^n$ &      &$E_1^n$  &           &$E_0^n$ &     &$E_1^n,F_{11}^n$ \\[1mm]
$3.3$ & $1$ & $0$ &  $1$ &$E_0^n,F_{00}^n$ & $E_0^n$ &$E_0^n$  & $E_0^n,F_{01}^n$  &$E_0^n$ &$E_0^n$ & $E_0^n,F_{00}^n$\\[1mm]
$3.4$ & $1$ & $1$ &  $1$ &$E_0^n,F_{00}^n$ &      &$E_1^n$  &           &$E_0^n$ &     &$E_1^n,F_{11}^n$  \\[1mm]
\hline
\end{tabular}
\end{center}
\label{Ta:333actions}
\end{table}}
\end{thm}

\begin{pf}The proofs of above three theorems are similar and therefore we give
the argument concerning the $(2,3,6)$-action only. By
\cite{[T]}, Theorem~4.1, such an action arises by the lifting an
automorphism group $\tilde G$ of an elliptic curve with a
presentation which differs from $\tilde G_C$ in
Table~\ref{Ta:reducedgroups} by having the additional relation
$\tilde x^m=\tilde y^{mk}$, for some integers $m, k$ such that $m$
divides $n$ and $k^2-k+1\equiv 0\;(n/m)$.    We shall prove that
the symmetric character of $X$ requires $m=n$ or $m=n/3$, which
leads to the presentations (\ref{C236}) or (\ref{D236}),
respectively. Here
$\theta(x_1)=c^3x,\theta(x_2)=c^2y,\theta(x_3)=c$ and it is easy
to check that the assignment $\theta(x_i)\rightarrow
\theta(x_j)^{-1},\theta(x_j)\rightarrow \theta(x_i)^{-1}$ does
not induce an automorphism of $G$ for any two distinct
$i,j\in\{1,2,3\}$. So assume that $\varphi :G\rightarrow G$ is
induced by the assignment
$\theta(x_1)\rightarrow\theta(x_1)^{-1}$ and
$\theta(x_2)\rightarrow\theta(x_2)^{-1}$. Then
$\varphi(x)=x^{-1},\varphi(y)=x^{-1}y\rho^\beta$ and
$\varphi(c)=c^{-1}y^{-1}x\rho^\alpha$. Thus
$\varphi(x)^m=\varphi(y)^{mk}\rho^\delta$ implies $(\tilde
x)^{-m}=(\tilde x^{-1}\tilde y)^{mk}=\tilde x^{-mk}\tilde
y^{mk}=\tilde x^{m-mk}$, where $\tilde x$ is the image of $x$ in
$\tilde G=G/\langle\rho\rangle$. Since $\tilde x^{2m}=\tilde
x^{mk}$, it follows that $k\equiv 2\;(n/m)$. Furthermore,
$k^2-k+1\equiv 0\;(n/m)$ which implies $3\equiv 0\;(n/m)$. Thus
$m=n$ or $m=n/3$. If $m=n$ then $G$ has the presentation
(\ref{C236}) and $\tilde G\cong\tilde G_C$. In the case $m=n/3$
it is convenient to replace the generators $x,y$ by $w=xy$ and
$v=y$. Since $x^{n/3}=y^{2n/3}\rho^\delta$, it follows that
$x^n=\rho^\delta$. On the other hand,
$x^n=cx^nc^{-1}=y^n\rho^{n\alpha}=\rho^{\varepsilon +n\alpha}$
which implies that $\rho^\delta =\rho^{\varepsilon +n\alpha}$.
Thus $w^{n/3}=x^{n/3}y^{n/3}\rho^{\gamma
n/3(n/3-1)/2}=\rho^{\delta +\varepsilon+ \gamma
n/3(n/3-1)/2}=\rho^{n\alpha +\gamma n/3(n/3-1)/2}$. Consequently
$G$ has the presentation (\ref{D236}) and $\tilde G\cong\tilde
G_D$.

The expression for the signature of $\Lambda$ in terms of the
parameters $\alpha,  \beta, \mu$  is an easy computation following
from $a_i = \varepsilon_i- 1$, $i=1,2,3$ and the relations
$\alpha=a_1$, $\beta=a_2$, $\mu = a_3$ and $\gamma
=a_1+a_3=\alpha +\mu$, which were proved in \cite{[T]}.
\end{pf}

\section{Full  conformal actions:  a corrigendum}\label{S:fullcorrigendum} 

  A 
Fuchsian group $\Lambda(\tau)$ with signature $\tau$  is 
 {\it finitely maximal\/} if  there is no Fuchsian group $\Lambda({\tau^\prime})$
with signature $\tau^\prime$ containing it with finite index.  Most Fuchsian groups are finitely maximal, but there a few infinite families of signatures, together with a few other individual signatures, for which the corresponding Fuchsian groups are not finitely maximal (see  \cite{[Sing2]}).   Suppose $\Lambda(\tau)$ is not finitely maximal.  
Let  $\Lambda(\tau^\prime)$ be a Fuchsian group such that 
$\Lambda(\tau) < 
\Lambda(\tau^\prime)$ with finite index.  Let $G <G^\prime$ be
finite groups such that $[G^\prime: G]=[\Lambda(\tau^\prime):
\Lambda(\tau)]$.  An
epimorphism
$\theta:
\Lambda(\tau)
\rightarrow G$ (with kernel  a surface group of genus $g$) determines a $G$-action on $X={\cal H}/\hbox{ker}(\theta)$.  $\theta$ may extend to an epimorphism
$\theta^\prime:
\Lambda(\tau^\prime)
\rightarrow G^\prime$ (with the same kernel).  If it does, 
$G$ is not the full automorphism group
of the surface $X= {\cal H}/\hbox{ker}(\theta)$.    When
$\tau, \tau^\prime$ are both triangular,  the $G$-action and the $G'$-action both determine the same exceptional point $X$.
  
 Clearly, if $G$ acts on an elliptic-hyperelliptic surface with a signature that is either finitely maximal or extends with finite index only to signatures which are not elliptic-hyperelliptic, it must be the full group of conformal automorphisms of the surface.      It turns out that the converse is also true:  no other triangular signature is the signature of the full group of conformal automorphisms of an elliptic-hyperelliptic surface.    This is not obvious, but is a  consequence of the following theorem, which is a corrigendum of Theorem~8.1 in \cite{[T]}.  
    
  \begin{thm}\label{T:corrigendum}
A group $G=\Lambda /\Gamma$ is the full group of conformal automorphisms of an
elliptic-hyperelliptic Riemann surface $X={\mathcal H}/\Gamma$ of
genus $g>5$ if and only if $G$ is one of the groups listed in the
Theorems ~3.1-7.1 of \cite{[T]} and $\Lambda$ has a finitely maximal
signature or $\Lambda$ has one of the following non-maximal signatures:
$[2,2,4,4],[4,4,8],[2,4,8],[2,2,8,8],[2,6,6],[2,2,6,6]$
corresponding to the ca\-ses $4.1,t=1;4.5,t=0;
4.4,t=0;4.6,t=1;6.2,t=0;6.2,t=1$, respectively, in \cite{[T]}.
\end{thm}

\section{Exceptional points}\label{S:exceptional}

In ${\cal M}_g^0$, $g>30$, there are exactly three exceptional points \cite{We04}, and  Singerman's results \cite{[Sing1]}  imply that these three exceptional points are also symmetric.  We show in this section  that the elliptic-hyperelliptic locus ${\cal M}_g^1$ can contain an arbitrarily large number of exceptional points,  but no more than four of them are also symmetric.      

  Let $G$ be the full group of conformal automorphisms of a surface $X \in {\cal M}_g^1$, and let $\rho$ denote the $1$-hyperelliptic involution.  Suppose $G/\langle \rho \rangle$ is of  the form (\ref{E:ellipticaut}) and $G$ acts with a triangular signature, so that $X$ is an exceptional point.  Then $G$ has order $2nmt$ and $t>2$.  
  Let ${\cal
R}_t$ denote the ring ${\cal G}$ if $t=4$ and ${\cal E}$ if
$t=3$ or $6$.  We recall that the prime factors of $n/m$ must split in ${\cal R}_t$.
\begin{lema}
Let  $p>3$ be a splitting prime factor of $n/m$ with multiplicity $\mu >0$ (i.e, $p^\mu$ is the highest power of $p$ dividing $n/m$).  Then the $p$-Sylow subgroup of (\ref{E:ellipticaut}), and hence of $G$, has  $1+\lfloor \mu/2 \rfloor$ distinct possible isomorphism types.
\end{lema}\label{L:splitSylow}
\begin{pf} Since $p$ splits in ${\cal R}_t$, $p=\pi \cdot \pi'$, and since $p>3$, $\pi$, $\pi'$ are non-associate primes which therefore generate distinct ideals $\pi{\cal R}_t$, $\pi'{\cal R}_t$.   Thus the $p$-Sylow subgroup of (\ref{E:ellipticaut}), which is also the $p$-Sylow subgroup of $G$, could be  any one of the following:
\begin{equation}
{{\cal 
R}_k\over \pi^l(\pi^\prime)^{\mu-l}{\cal R}_k}\simeq {{\cal 
R}_k\over \pi^l{\cal R}_k} \oplus{{\cal 
R}_k\over (\pi^\prime)^{\mu-l}{\cal R}_k}\simeq 
{\Bbb Z}_{p^{l}} 
\oplus {\Bbb Z}_{p^{\mu-l}},\quad 0\leq l
\leq \lfloor
\mu/2\rfloor.
\end{equation}
\end{pf}

\begin{thm}\label{T:infsequence}
There exist infinite sequences of genera in which the number of exceptional points in ${\cal M}_g^1$ is larger than any pre-assigned positive integer.
\end{thm}

\begin{pf}
If $G$ is the full group of conformal automorphisms of an exceptional point $X \in {\cal M}_g^1$, it acts with one of signatures $[2,4,8]$, $[4,4,8]$,
$[2,6,6]$, $[2,3,12]$, $[4,3,6]$, $[4,6,6]$ or $[4,6,12]$,  and by
(\ref{E:RiemHur}), the ratio $[nm : (g-1)]$ is
 $$[a: b] \in \{[2:1], [2:3], [1:1],[2:3], [2:5],[1:3]\}.$$
   These ratios place restrictions on the multiplicities of prime divisors of $g-1$.  If $p$ is a nonsplitting prime in ${\cal R}_t$, the  $p$-Sylow subgroup of ${\Bbb Z}_n \oplus {\Bbb Z}_m$ is  of the form ${\Bbb Z}_{p^k} \oplus {\Bbb Z}_{p^k}$ by Theorem~\ref{T:pSylow}.  Thus if $p>\hbox{max}\{a,b\}$, it must have even (possibly $0$) multiplicity in $g-1$.    There is no restriction on the multiplicity of a splitting prime $p>\hbox{max}\{a,b\}$ in $g-1$.   If the multiplicity of such a splitting prime is $\mu$,
  by Lemma~\ref{L:splitSylow}, there are at least $1 + \lfloor \mu/2 \rfloor$ different isomorphism types for $G$, hence at least that many distinct exceptional points in ${\cal M}_g^1$.

Given $a,b$, $t$, an ${\cal R}_t$-splitting prime $p>5$, and an arbitrary positive integer $N$, there exist infinite sequences of genera $g$ such that  the multiplicity of $p$   in $g-1$  is greater than N.  It follows that there are infinitely many genera $g$ in which the number of exceptional points in the elliptic-hyperelliptic locus is larger than $N$.
\end{pf}

It is easy to construct infinite sequences of genera in which the elliptic-hyperelliptic locus contains no exceptional points.  For example,  $g_n=1+p^{2n+1}$, $n=0,1, \dots$, where $p$ is any prime such that $p \equiv -1$ mod $12$, is such a sequence, since $p$ has odd multiplicity in $g_n-1$,  and does not split in ${\cal G}$ or in ${\cal E}$.  

We now restrict attention to exceptional points which are also symmetric.

\begin{thm}\label{T:genus}
The genus $g$ of a symmetric exceptional point $X \in {\cal M}_g^1$  is $ka^2+1$ for some integer $a$ and
$k\in\{1,2,3,6,10,30\}$. For such $g$ there are the following
non-equivalent triangular symmetric actions on $X$; those  in bold type correspond to the full automorphism groups.
\begin{enumerate}
\item $k=1$ and $a\equiv 0\;(2):$ $4.2.B_{00}^a,B_{01}^a,$
$4.3.A_0^a,A_1^a,$ ${\bf 4.4.B_{00}^{2a},B_{01}^{2a},}$
$4.6.A_{0}^a,A_{1}^a,$ ${\bf 6.2.C_0^a},\;\;3.3.E_0^a$; and $ {\bf
6.8.D_{00}^a}$ if $a\equiv 0\;(3)$,
\item $k=1$ and $a\equiv 1\;(2):$
${\bf 6.2.C_1^a}$, $3.3.E_0^a$; and ${\bf 6.8.D_{01}^a}$ if
$a\equiv 0\;(3)$,
\item $k=2$ and $a\equiv 0\;(2):$
$4.2.A_0^a,A_1^a$, $4.3.B_{00}^{2a},$ $B_{01}^{2a},$ ${\bf
4.4.A_0^{2a},A_1^{2a}},$ $4.6.B_{00}^{2a},B_{01}^{2a},$ ${\bf
6.3.C_0^{2a}},$ $6.7.C_0^a,$ $3,2.E_0^{2a}$;  and
$6.4.D_{00}^{2a},\;{\bf 6.5.D_{00}^{2a}},\;3.4.F_{00}^{2a}$ if
$a\equiv 0\;(3)$,
\item $k=2$ and $a\equiv 1\;(2):$
$4.2.A_0^a,A_1^a,\;4.3.B_{11}^{2a},\;B_{10}^{2a},$ ${\bf
4.4.A_0^{2a},A_1^{2a}},$ $4.6.B_{10}^{2a},$ $B_{11}^{2a},$
$6.7.C_1^a$ ${\bf 6.3.C_1^{2a}}$ $3,2.E_1^{2a}$;  and
$6.4.D_{11}^{2a},$ ${\bf 6.5.D_{11}^{2a}},$ $ 3.4.F_{11}^{2a}$ if
$a\equiv 0\;(3)$,
\item $k=3$ and $a\equiv 0\;(2):$
${\bf 4.5.B_{00}^{2a}},$ ${\bf B_{01}^{2a}},$ ${\bf 6.8.C_0^a},$
${\bf 6.2.D_{00}^{3a}},$ $3.3.F_{00}^{3a}$,
\item $k=3$ and $a\equiv 1\;(2):$
${\bf 6.8.C_0^a,}$ ${\bf 6.2.D_{10}^{3a}},$ $3.3.F_{01}^{3a},$
\item $k=6$ and $a\equiv 0\;(2):$
${\bf 4.5.A_0^{2a},A_1^{2a}},$ ${\bf 6.3.D_{00}^{6a}},$ $
6.4.C_0^{2a}$, ${\bf 6.5.C_0^{2a}},$ $6.7.D_{00}^{3a},$
$3.2.F_{00}^{6a},$ $3.4.E_0^{2a},$
\item $k=6$ and $a\equiv 1\;(2):$
${\bf 4.5.A_0^{2a},A_1^{2a}}$, ${\bf 6.3.D_{11}^{6a}},$
$6.4.C_1^{2a}$, ${\bf 6.5.C_1^{2a}},$ $6.7.D_{11}^{3a},$
$3.2.F_{11}^{6a},$ $3.4.E_1^{2a},$
\item $k=10:$ ${\bf 6.6.C_0^{2a}}$ or ${\bf 6.6.C_1^{2a}}$
according to $a$ being even or odd,
\item $k=30:$ ${\bf 6.6.D_{00}^{6a}}$ or
${\bf 6.6.D_{11}^{6a}}$ according to $a$ being even or odd.
\end{enumerate}
\end{thm}
\begin{pf}
Let $G=\Lambda/\Gamma$ be a group of conformal automorphisms acting with a triangular signature on a
symmetric elliptic-hyperelliptic  surface $X= {\mathcal
H}/\Gamma$ of genus $g$.  Thus $\Gamma$ is a surface group of genus $g$ and $\Lambda$ has one of the signatures (\ref{E:ehypsignatures}).  Let $\rho$ be the $1$-hyperelliptic
involution of $X$. Theorems~\ref{T:action244}, \ref{T:action236}, and \ref{T:action333} show that the reduced group $G/\langle \rho \rangle$ is isomorphic to one of groups listed
in Table~\ref{Ta:reducedgroups}.  The genus of $X$ is determined as in the following example.  Suppose  $G/\langle \rho \rangle\cong\tilde
G_D$.  Then the order of $G$ is $12n^2/3=4n^2$ and the Riemann-Hurwitz relation (\ref{E:RiemHur}) where $\Lambda$ has  signature  $[2\varepsilon_1, 3\varepsilon_2, 6 \varepsilon_3]$ with at least one $\varepsilon_i$ equal to $2$, is $g= 1 + 2n^2(1 - {1}/{2\varepsilon_1} - {1}/{3\varepsilon_2} - {1}/{6\varepsilon_3})$.  The expression in parentheses is either $1/4, 5/12, 1/2, 1/6, 1/4$ or $1/12$.  Since $g$ is an integer, $n$ must be divisible by $2$, $6$, $1$, $3$, $2$, $6$, respectively.  Thus there exists an integer $a$ such that $g=ka^2 + 1$, for $k=2, 30, 1, 3, 6$, respectively.   Similar computations using the other reduced groups $G/\langle \rho \rangle$ in Table~\ref{Ta:reducedgroups}, and compatible signatures, yield the general result $g=ka^2+1$ for some  
$k\in\{1,2,3,6,10,30\}$.

The actions are determined as follows:  Given a signature and a reduced group $G/\langle \rho \rangle$  (compatible with the signature),  we use the tables in Theorems~\ref{T:action244}, \ref{T:action236}, and \ref{T:action333},  together with the value of $k$ calculated as above.      For example,  actions with signature $[4,6,6]$ have reduced group  $\tilde G_D$ or $\tilde G_C$, and correspond to Case 6.6  in Table~\ref{Ta:236actions}.  The Riemann-Hurwitz relation using $\tilde G_D$ implies   $n$ is divisible by $6$ and $k=30$.  Putting $n=6a$, we have $n \equiv 0$ or $\equiv 6$ mod $12$ according as $a$ is even or odd.  Thus the actions are $6.6 D^{6a}_{00}$ if $a$ is even and $6.6. D^{6a}_{11}$ if $a$ is odd.  Similarly, the Riemann-Hurwitz relation using $\tilde G_C$ implies $n$ is divisible by $2$ and $k=10$.  Thus $n \equiv 0, \pm 2, \pm 4$ or $6$ mod $12$.  From Table~\ref{Ta:236actions} the actions are either $6.6.C^{2a}_0$ or $6.6.C^{2a}_1$.    Since $[4,6,6]$ does not extend to an elliptic-hyperelliptic signature,  these actions are  full actions.   This yields statements 9 and 10 of the theorem.  The other statements are derived similarly.  We appeal to Theorem~\ref{T:corrigendum} to determine  extensions of $G$-actions with signature $\tau$ to $G'$-actions with signature $\tau'$.  These are given in Table~\ref{Ta:fullextensions},  below. 
\end{pf}

\begin{coro} If ${\cal M}_g^1$ contains exceptional points, at most four of them are symmetric.
\end{coro}
\begin{pf} Full actions are in bijection with exceptional points.  
 If $g-1=ka^2=k'a'^2$, for $k,k' \in \{1,2,3,6,10, 30\}$, then  $k=k'$ and hence also $a=a'$.  Thus the ten enumerated cases in Theorem~\ref{T:genus} are mutually exclusive.  One merely observes that the maximum number of full actions in any one case is four.
\end{pf}

\begin{table}[htdp]
\caption{Extensions of triangular actions on symmetric surfaces in ${\cal M}_g^1$}
\begin{center}
\begin{tabular}[!htp]{||l|l|l|l|l|c|}
\hline
&&&&&\\[-3mm]
$g$ & $\tau$ &$\;\;\;\;\;\;\;G$ &$\tau'$ &$\;\;\;\;\;\;\;G'$& $\tau':\tau$\\[-3mm]
&&&&&\\
\hline
&&&&&\\[-3mm]
 $a^2+1$&$[4,8,8]$ & $4.2.B_{00}^a,B_{01}^a$&$[2,8,8]$&
$4.6.A_0^a,A_1^a$&$2$\\[1mm]
 $2a^2+1$& $[4,8,8]$ & $4.2.A_0^a,A_1^a$& $[2,8,4]$ &
$4.4.A_0^{2a},A_1^{2a}$& $4$\\[1mm]
 $2a^2+1$&$[4,8,8]$& $4.2.A_0^a,A_1^a$ & $[2,8,8]$ &
$4.6.B_{00}^{2a},B_{01}^{2a}$ & $2$\\[1mm]
$a^2+1$ & $[4,4,4]$&$4.3.A_0^a,A_1^a$& $[2,4,8]$ & $4.4.B_{00}^{2a},B_{01}^{2a}$ & $2$\\[1mm]
 $2a^2+1$ & $[4,4,4]$ & $4.3.B_{00}^{2a},B_{01}^{2a}$&
$[2,4,8]$ & $4.4.A_0^{2a},A_1^{2a}$& $2$\\[1mm]
 $a^2+1$ & $[2,8,8]$ & $4.6.A_0^a,A_1^a$& $[2,8,4]$
&$4.4.B_{00}^{2a},B_{01}^{2a}$&$2$\\[1mm]
 $2a^2+1$& $[2,8,8]$ & $4.6.B_{00}^{2a},B_{01}^{2a}$&
$[2,8,4]$ & $4.4.A_0^{2a},A_1^{2a}$& $2$\\[1mm]
 $2a^2+1$ & $[2,6,12]$ & $6.4.D_{00}^{2a}$ & $[2,3,12]$ &
$6.3.C_0^{2a}$ & $3$\\[1mm]
 $6a^2+1$ & $[2,6,12]$ & $6.4.C_0^{2a}$ & $[2,3,12]$ &
$6.3.D_{00}^{6a}$ & $3$\\[1mm]
 $6a^2+1$ & $[4,3,12]$ & $6.7.D_{00}^{3a}$ & $[2,3,12]$ &
$6.3.D_{00}^{6a}$ & $4$\\[1mm]
 $2a^2+1$ & $[4,3,12]$ & $6.7.C_0^a$ & $[2,3,12]$ &
$6.3.C_0^{2a}$ & $4$\\[1mm]
 $2a^2+1$ & $[6,3,3]$ & $3.2.E_0^{2a}$ & $[2,3,12]$ &
$6.3.C_0^{2a}$ & $2$\\[1mm]
 $6a^2+1$ & $[6,3,3]$ & $3.2.F_{00}^{6a}$ & $[2,3,12]$ &
$6.3.D_{00}^{6a}$ & $2$\\[1mm]
 $3a^2+1$&$[3,6,6]$ & $3.3.F_{00}^{3a}$ & $[2,6,6]$ &
$6.2.D_{00}^{3a}$ & $2$\\[1mm]
 $a^2+1$&$[3,6,6]$ & $3.3.E_0^a$ & $[2,6,6]$ & $6.2.C_0^a$ &$2$\\[1mm]
  $2a^2+1$&$[6,6,6]$ & $3.4.F_{00}^{2a}$& $[6,3,3]$ &
$3.2E_0^{2a}$& $3$\\[1mm]
\hline
\end{tabular}
\end{center}
\label{Ta:fullextensions}
\end{table}

\section{Symmetries of exceptional points}\label{S:symmetries}

 Let
$X={\mathcal H}/\Gamma$ be a symmetric elliptic-hyperelliptic
Riemann surface whose full group of conformal automorphisms
$G=\Lambda/\Gamma$ acts with a triangular signature $[k,l,m]$. The existence of a symmetry on $X$ means that $\Lambda$ is the canonical Fuchsian group of a proper 
 NEC group $\widetilde \Lambda$, containing $\Lambda$ with index $2$, and containing $\Gamma$ as a normal subgroup.   Then $\tilde\Lambda/\Gamma = {\cal A}$ is the full group of conformal and anticonformal automorphisms of $X$.   Comparison of (\ref{sig}) and (\ref{sig-fuch}) yields two possibilities for the signature of $\widetilde \Lambda$:  
$(0;+;[-];(k,l,m))$, and, if $k=l$, 
$(0;+;[k];\{(m)\})$; we shall see that, in this context, the second possibility does not occur.
   
Let  $\tilde \theta$ be the canonical
epimorphism $\widetilde \Lambda\rightarrow \mathcal{A}$ with kernel $\Gamma$.  A symmetry $\phi \in {\mathcal A}$ is the image
under $\widetilde \theta$ of an element $d$ from the subset $\widetilde
\Lambda\setminus\Lambda$ of orientation-reversing elements of $\widetilde \Lambda$. If $d$ cannot be chosen as a reflection
then $\phi$ has no ovals.  Otherwise $d$ is conjugate
to one of the reflection generators $c$ in the canonical system of generators of $ \widetilde \Lambda$.  The
number of ovals $||\phi||$ is the number of empty period
cycles in the group $\widetilde \Gamma=\widetilde \theta^{-1} (
\langle \phi\rangle )$.    A formula for $||\phi||$ is given in \cite{[G1]} in terms of orders of centralizers:  
\begin{equation}\label{ovals}
||\phi||= \sum {|C( {\mathcal A},\widetilde  \theta(c_i))|
/|\widetilde \theta(C(\widetilde  \Lambda,c_i))|,}
\end{equation}
where $c_i$ runs over  pairwise non-conjugate canonical
reflection generators in $\widetilde  \Lambda$ whose images are conjugate to
$\phi$, and $C(A,a)$ denotes the centralizer of the element $a$ in the group $A$.   In \cite{[Sing3]} (see also \cite{[Sing1]}) it is proved that the centralizer of a reflection $c$ in an NEC group $\Lambda$ is  isomorphic to ${\Bbb Z}_2 \oplus {\Bbb Z}$ if the associated period cycle in $\Lambda$ is empty or consists of odd periods only; otherwise it is isomorphic to ${\Bbb Z}_2 \oplus({\Bbb Z} * {\Bbb Z})$, where $*$ denotes the free product.

Using  these results we obtain the following classification of the centralizers of reflections in an NEC group whose canonical Fuchsian group is a triangle group.  The notation $c_1 \sim c_2$ denotes conjugacy in $\widetilde \Lambda$.

\begin{lema}\label{sing1}
(a) Let $\widetilde  \Lambda$ be an NEC group with signature
$(0;+;[-];(k',l',m'))$ and let $c_0,c_1,c_2$ be the canonical system of
generators of $\widetilde  \Lambda$. Then
 
\medskip
\noindent $(i)$ for $k'=2k+1,l'=2l+1,m'=2m+1,  c_0\sim c_1\sim
c_2$ and

$\; C(\widetilde  \Lambda,c_0)= \langle c_0 \rangle \oplus(
\langle (c_2c_0)^{m}(c_1c_2)^{l}(c_0c_1)^{k} \rangle )$,

\medskip
\noindent $(ii)$ for $k'=2k,l'=2l+1,m'=2m+1, c_0\sim c_1\sim c_2$
and

$\; C(\widetilde  \Lambda,c_0)= \langle c_0 \rangle \oplus (
\langle (c_0c_1)^{k} \rangle \ast \langle (c_2c_0)^{m}(c_1c_2)^{l}
(c_1c_0)^{k}(c_2c_1)^{l}(c_0c_2)^{m} \rangle  )$,

\medskip
\noindent $(iii)$ for $k'=2k,l'=2l,m'=2m+1, c_0\sim  c_2$ and

$\; C(\widetilde  \Lambda,c_0)= \langle c_0 \rangle \oplus(
\langle (c_0c_1)^{k} \rangle \ast \langle
(c_2c_0)^{m}(c_2c_1)^{l}(c_0c_2)^{m} \rangle )$,

$\; C(\widetilde  \Lambda,c_1)= \langle c_1 \rangle \oplus(
\langle (c_0c_1)^{k} \rangle \ast \langle (c_1c_2)^{l} \rangle )$,

\medskip \noindent
$(iv)$ for $k'=2k,l'=2l,m'=2m$

$\; C(\widetilde  \Lambda,c_0)= \langle c_0 \rangle \oplus(
\langle (c_0c_1)^{k} \rangle \ast \langle (c_0c_2)^{m} \rangle )$,

$\; C(\widetilde  \Lambda,c_1)= \langle c_1 \rangle \oplus(
\langle (c_0c_1)^{k} \rangle \ast \langle (c_1c_2)^{l} \rangle )$,
$\; C(\widetilde  \Lambda,c_2)= \langle c_2 \rangle \oplus(
\langle (c_0c_2)^{m} \rangle \ast \langle (c_1c_2)^{l } \rangle
)$.

\medskip \noindent

(b) Let $\widetilde  \Lambda$ be an NEC group with signature
$(0;+;[k];\{(m)\})$ and let $x,e,c_{0},c_{1}$ be a canonical
system of generators of $\widetilde \Lambda$.  Then $c_0 \sim c_1$ and 
$$\;\;
C(\widetilde  \Lambda,c_{0})= \left\{
\begin{array}{lr}
<c_{0}>\oplus <(c_{0}c_{1})^{m/2}>\ast<e(c_{0}c_{1})^{m/2}e^{-1}>
& if \;\;m \;is\; even ,\\
<c_{0}>\oplus(<e(c_{0}c_{1})^{(m-1)/2}>) & if \;\;\;\; m \;is
\;odd .
\\
\end {array}
\right. $$
\end{lema}

\begin{thm}\label{T:fullgroups}
Let $\mathcal{A}$ be the full group of conformal and anticonformal automorphisms of a symmetric exceptional point $X \in {\cal M}_g^1$ with 
genus  $g=ka^2+1$,  $k\in \{1,2,3,6,10,30\}$.  Then ${\mathcal A}$ acts with NEC signature $(0;+;[-];(k',l',m'))$ and is a semidirect product
$G \rtimes \langle\tau :\tau ^2\rangle$, where  $G$ is  
the full group of conformal automorphisms of $X$, listed in Table~\ref{Ta:full} according to the value of $k$, $\tau$ is a symmetry of $X$, and the action of $\tau$ on the generators of $G$ is given in 
Table~\ref{Ta:tauaction}.  
$\mathcal{A}$ contains two or three conjugacy classes of symmetries with fixed points, the number of whose ovals is 
given in Table~\ref{Ta:full}.

\begin{table}[htdp]
\caption{Conformal group $G$ and ovals}
\begin{center}
\begin{tabular}[!htp]{|l|l|l|l||l|l|l|l|} \hline
$k$ & $\;\;\;\;G$ & $\;\;$Ovals & Conditions& $k$ & $\;\;\;\;G$ & $\;\;\;$Ovals & Conditions\\
\hline
&&&&&&&\\[-3mm]
$1$ & $4.4.B_{00}^{2a}$ & $a,2,a$& $a\equiv 0\;(2)$ &$3$ &
$4.5.B_{01}^{2a}$&$2a,a,a$ & $a\equiv
0\;(2)$\\[1mm]
 & $4.4.B_{01}^{2a}$ & $a,2,a$&
$a\equiv 0\;(2)$& & $4.5.B_{00}^{2a}$&$2a,a,a$ & $a\equiv
0\;(2)$\\[1mm]
 & $6.2.C_0^a$ & $2,a,a$&
$a\equiv 0\;(2)$&&$6.2.D_{00}^{3a}$ & $2,a,a$&$a\equiv
0\;(2)$\\[1mm]
 & $6.2.C_1^a$ & $1,a,a$&
$a\equiv 1\;(2)$&&$6.2.D_{10}^{3a}$ & $1,a,a$&$a\equiv
1\;(2)$\\[1mm]
 & $6.8.D_{00}^a$ & $a,\frac{2}{3}a$&
$a\equiv 0\;(6)$&&$6.8.C_0^a$ & $a,2a$&none\\[1mm]
 & $6.8.D_{01}^a$ & $a,\frac{2}{3}a$&
$a\equiv 3\;(6)$&$6$ & $4.5.A_0^{2a}$ & $4a,2a,a$ &
none\\[1mm]
 $2$ & $4.4.A_0^{2a}$ & $2a,4,a$ &
none&& $4.5.A_1^{2a}$ & $4a,2a,a$ &
none\\[1mm]
 & $4.4.A_1^{2a}$ & $2a,1,a$ &
none&&$6.3.D_{00}^{6a}$ & $a,3a$&$a\equiv
0\;(2)$\\[1mm]
 & $6.3.C_\varepsilon^{2a}$ & $a,a$ &$\varepsilon\equiv a\;(2)$
&&$6.3.D_{11}^{6a}$ & $a,3a$&$a\equiv
1\;(2)$\\[1mm]
& $6.5.D_{00}^{2a}$ & $\frac{1}{3}a,a$ & $a\equiv 0\;(6)$&&
$6.5.C_\varepsilon^{2a}$ & $a,a$ & $\varepsilon\equiv
a\;(2)$\\[1mm]
& $6.5.D_{11}^{2a}$ & $\frac{1}{3}a,a$ & $a\equiv 3\;(6)$&$30$ &
$6.6.D_{00}^{6a}$ & $3a,3a$ & $a\equiv
0\;(2)$\\[1mm]
$10$ & $6.6.C_\varepsilon^{2a}$ & $a,3a$ & $\varepsilon\equiv
a\;(2)$&& $6.6.D_{11}^{6a}$ & $3a,3a$ & $a\equiv
1\;(2)$\\[1mm]
\hline
\end{tabular}
\end{center}
\label{Ta:full}
\end{table}

\begin{table}[htdp]
\caption{Action of $\tau$ on $G$}
\begin{center}
\begin{tabular}[!htp]{|l|l|l|l|l|l|l|} \hline
$\;\;G$ & $\;\;\;\;\;\;\;\tau c\tau$ & $\tau x\tau$ & $\tau y\tau$
& $\;\;\;\;\;\;\tau w\tau$ & $\;\;\;\;\;\;\;\;\tau v\tau$ & $\tau\rho\tau$ \\
\hline
&&&&&&\\[-3mm]
$A_{\varepsilon}$ & $c^{-1}$ & $x\rho^{\alpha}$& $y^{-1}$ & & & $\rho$\\[1mm]
$B_{\varepsilon\delta}$ & $c^{-1}$ & & &$v^{-2}w^{-1}\rho^{\alpha +\gamma}$&$v\rho^\alpha$ & $\rho$ \\[1mm]
$C_\varepsilon$ & $y^{-1}c^{-1}\rho^{\alpha +\beta}$ & $x^{-1}$ &
$cyc^{-1}$ & & & $\rho$\\[1mm]
$D_{\varepsilon\varepsilon'}$ & $v^{-1}c^{-1}\rho^{\alpha+\beta}$
& & &$v^3w^{-2}\rho^{\beta}$ &
$v^2w^{-1}\rho^{\alpha+\beta+\mu}$& $\rho$\\[1mm]
\hline
\end{tabular}
\end{center}
\label{Ta:tauaction}
\end{table}

\end{thm}

\begin{pf}  Let $ {\mathcal A}=\tilde \Lambda /\Gamma$ and let $\tilde \theta:\widetilde \Lambda\rightarrow \mathcal{A}$ be the canonical epimorphism with kernel $\Gamma$.  We first show that $\tilde \Lambda$ cannot have the signature $(0;+;[k];\{(m)\})$.  Suppose it had such a signature, and let $x,e, c_0, c_1$ be a canonical system of generators.  Then $x_1=x,x_2=c_0x^{-1}c_0,x_3=c_0c_1$ is a system of canonical generators of the canonical Fuchsian group $\Lambda$, which has signature $[k,k,m]$.    The full group of conformal automorphisms $\Lambda/\Gamma =G$ is generated by  $g_1=\tilde
\theta(x_1), g_2=\tilde \theta(x_2)$, both of order $k$.   Since
$c_0x_1c_0=x_2^{-1}$ and $c_0x_2c_0^{-1}=x_1^{-1}$, it follows, for $\tau = \tilde\theta(c_0)$, 
that
 \begin{equation}\label{second type}\tau
g_1\tau=g_2^{-1}\;\;{\rm and }\;\; \tau
g_2\tau=g_1^{-1}.\end{equation}
 Thus $g_1 \mapsto g_2^{-1}$, $g_2 \mapsto g_1^{-1}$ is an 
outer automorphism of $G$.   This is false if $G$ has type 4.5, 6.2 or 6.6.  On the other hand,  these are the only cases where the action is full and the signature has  the form $[k,k,m]$.

Thus we may assume $\tilde
\Lambda$ has signature  $(0;+;[-];(k',l',m'))$.  Let 
$c_0,c_1,c_2$ be the system of canonical generators of $\tilde
\Lambda$. Then $x_1=c_0c_1,x_2=c_1c_2$ and $x_3=c_2c_0$ is a system of canonical generators of the canonical Fuchsian group $\Lambda$. Since
$c_1x_1c_1=x_1^{-1}$ and $c_1x_2c_1=x_2^{-1}$, it follows that
\begin{equation}\label{first type}\tau g_1\tau=g_1^{-1}\;\;{\rm and}\;\;\tau
g_2\tau=g_2^{-1},\end{equation}
for $g_1=\tilde
\theta(x_1),g_2=\tilde \theta(x_2)$ and $\tau =\tilde
\theta(c_1)$. Thus $\mathcal{A}$ is a semidirect product
$G \rtimes {\Bbb Z}_2=\langle g_1,g_2\rangle \rtimes \langle\tau
:\tau^2\rangle$, where $\tau$ is a symmetry of $X$. $G$ is the full group of conformal automorphisms of $X$, hence it must have one of the actions 
listed in bold  in Theorem~\ref{T:genus}.  $G$ has generators $x,y,c$ such that
$g_1,g_2,(g_1g_2)^{-1}$ are equal to $c^{-2}x,c,y^{-1}c$, respectively,  if $G$ has a $(2,4,4)$-action, or to 
$c^3x,c^2y,c$, if $G$ has a  $(2,3,6)$-action.
(\ref{first type}) induces the action of $\tau$ on
$x,y,c$ given in Table~\ref{Ta:tauaction}.

A symmetry of $X$ with fixed points is conjugate to one of
  $\tilde \theta(c_0)=g_1\tau$, $\tilde \theta(c_1)=\tau$ or
$\tilde \theta(c_2)=\tau g_2$ and we shall calculate the number of
its ovals using (\ref{ovals}). Let $\upsilon_i$ and
$\tilde \upsilon_i$ denote the orders of $\tilde \theta (C(\tilde
\Lambda ,c_i))$ and $C( {\mathcal A},\widetilde \theta(c_i))$,
respectively. Any element $h\in \mathcal{A}$ has the unique
presentation
\begin{equation}\label{element}
h=x^ry^sc^t\tau^p\rho^q,\end{equation} determined by a sequence
$(r,s,t,p,q)$ of integers satisfying $0\leq p,q \leq 1$, $0\leq t<
4$ (or $6$), $0\leq s<n$ and $0\leq r< d$, where $d=n$ for
$A_{\varepsilon}^n$ and $C_{\varepsilon}^n;$ $d=n/2$ for
$B_{\varepsilon\delta}^n$;  and $d=n/3$ for
$D_{\varepsilon\varepsilon'}^n$. Furthermore, for any integer
$\kappa$ \begin{equation}\label{translations}
(x^ry^s)^\kappa=x^{r\kappa}y^{s\kappa}\rho^{\gamma
\kappa(\kappa-1)/2}.
\end{equation}

Assume $G$ has $(2,4,4)$-action. Then we check that
any element of $C( {\mathcal A},g_1\tau)$ is determined by
$$ (r,s,0,p,q), \;{\rm with }\;x^{2r}=\left\{
\begin{array}{ll}
\rho^{s(\alpha +\gamma)} & \hbox{if\ }p=0,\\
\rho^{\alpha(s+1)+\mu+\gamma s} & \hbox{if\ }p=1;\\
\end {array}
\right.$$ 
or
 $$ (r,s,2,p,q), \;{\rm with}\;x^{2(r+1)}=\left\{
\begin{array}{ll}
\rho^{s(\alpha +\gamma)} & \hbox{if\ }p=1,\\
\rho^{\alpha(s+1)+\mu+\gamma s} & \hbox{if\ }p=0.\\
\end {array}
\right.$$ 
If $G=A_{\varepsilon}^n$ with even $n$, we
obtain the following possible sequences: \begin{enumerate}
\item
$(0,s,0,0,q),(-1,s,2,1,q):s(\gamma +\alpha)\equiv 0\;(2),$
\item
$(0,s,0,1,q),(-1,s,2,0,q):s(\gamma +\alpha)+\alpha+\mu\equiv
0\;(2)$
 \item $(n/2,s,0,0,q),(\frac{n}{2}-1,s,2,1,q):s(\gamma
+\alpha)\equiv \varepsilon\;(2),$
\item
$(n/2,s,0,1,q),(\frac{n}{2}-1,s,2,0,q):s(\gamma
+\alpha)+\alpha+\mu\equiv \varepsilon\;(2)$\end{enumerate} 
whose total 
number, given particular values of $\alpha,\gamma,\mu$ and
$\varepsilon$, is $\tilde \upsilon_0=8n$.  For
$G=A_{\varepsilon}^n$ with $n$ odd and
also for $G=B_{\varepsilon\delta}^n$, we have only the first two possibilities, so that 
 $\tilde \upsilon_0=4n.$

$C( {\mathcal A},\tau)=\{x^ry^s\tau^p\rho^q:
y^{2s}=\rho^{r\alpha}\}\cup
\{x^ry^sc^2\tau^p\rho^q:y^{2s}=\rho^{r\alpha +\mu}\}$, so its
elements are determined by the sequences
\begin{enumerate}\item $(r,0,0,p,q): \alpha r\equiv 0\;(2);$
$(r,\frac{n}{2},0,p,q):\alpha r\equiv \varepsilon\;(2);$ \item
$(r,0,2,p,q): \alpha r+\mu\equiv 0\;(2);$
$(r,\frac{n}{2},2,p,q):\alpha r+\mu \equiv \varepsilon\;(2).$
\end{enumerate} In all cases
but $4.4$, $\tilde \upsilon_1=\tilde \upsilon_0$ and in the
exceptional case $\tilde \upsilon_1=16n$ or $8n$ according to $G$
being $A_{0}^n$ or $A_{1}^n,B_{0\delta}^n$, respectively.

$C( {\mathcal A},\tau
g_2)=\{x^ry^s\rho^q,x^ry^sc^3\tau\rho^q:x^{r-s}y^{s-r}=\rho^{s\gamma+\alpha(r+s)}\}\cup
\{x^ry^sc^2\rho^q,$ $x^ry^sc\tau\rho^q:$
$x^{r-s}y^{s-r}=\rho^{s\gamma+\mu+\alpha(r+s)}\}$. Thus for
$G=A_{\varepsilon}^n$, $\tilde v_2$ is the total number of sequences of the form
$(r,r,0,0,q),$ $(r,r,3,1,q)$, $r\gamma\equiv 0\;(2)$ and $(r,r,2,0,q),(r,r,1,1,q)$,  $r\gamma+\mu \equiv 0\;(2)$.   According to the particular values of $\gamma$ and
$\mu$, $\tilde v_2=4n$, except the case $4.3$, where $\tilde
\upsilon_2=8n$. For $G=B_{\varepsilon\delta}^n$, only half
of the  listed sequences are possible; however we must now include sequences of the form 
$(r,\frac{n}{2}+r,0,0,q),(r,\frac{n}{2}+r,3,1,q)$ with $\alpha
n/2 +\gamma r+\delta \equiv 0\;(2)$ and
$(r,\frac{n}{2}+r,2,0,q),(r,\frac{n}{2}+r,1,1,q)$ with $\alpha
n/2+ r\gamma+\mu+\delta \equiv 0\;(2)$.  In all cases, except $4.3$,
there are $2n$ additional sequences and in the exceptional
case, $4n$ or $0$ according to $\varepsilon +\delta$ being even
or odd.  Consequently, $\tilde \nu=4n$, except the case $4.3$ with
$\varepsilon +\delta\equiv 0\;(2)$, where $\tilde \nu=8n$.

Let $k,l,m$ be integers such that $k'=2k,l'=2l$ and $m'=2m$. Then
by Lemma \ref{sing1},\label{centralizator}
\begin{equation}\label{all even}\upsilon_0=4\cdot\hbox{ord}(g_1^k(g_1g_2)^m),\upsilon_1=4\cdot\hbox{ord}
(g_1^kg_2^l)\;{\rm and
}\;\upsilon_2=4\cdot\hbox{ord}((g_1g_2)^mg_2^l).\end{equation} 
Thus we
obtain the following values of $(\nu_0,\nu_1,\nu_2)$:
$4.2\;(4,4,4)$; $4.3\;(8,8,2n\cdot2^{\varepsilon+\delta\;(2)})$
for $G=B_{\varepsilon\delta}^n$ and $(8,8,4n)$ for
$G=A_{\varepsilon}^n$; $4.4\;(8,4n\cdot 2^{\varepsilon},8)$;
$4.5\;(4,8,8)$ and $4.6\;(8,8,4)$.

Now we check that $h\tau g_2h^{-1}\neq g_1\tau$ and $h\tau
g_2h^{-1}\neq \tau$ for any element $h$ of the form
(\ref{element}) while $\tau$ and $g_1\tau$ are conjugate only for
odd $n$ via $x^{(n-1)/2}c^3$ or $x^{(n-1)/2}c$ according to
$\varepsilon=0$ or $\varepsilon =1$. Thus by (\ref{ovals}), there
are two or three conjugacy classes of symmetries with fixed points
and the numbers of their ovals are equal to $\tilde
\upsilon_0/\upsilon_0,\tilde \upsilon_1/\upsilon_1,\tilde
\upsilon_2/\upsilon_2$ or $\tilde \upsilon_0/\upsilon_0+\tilde
\upsilon_1/\upsilon_1,\tilde \upsilon_2/\upsilon_2$ according to
$n$ being even or odd.

The arguments are similar in the case that $G$ has a $(2,3,6)$ action, and we omit them.

\end{pf}

\end{document}